\documentclass[10pt,a4paper,draft]{article}
\usepackage{amsmath,amssymb,amsfonts, amsthm}
\date{\begin{flushleft} {\small } \end{flushleft}}
\usepackage[cp1251]{inputenc}
\usepackage[english,russian]{babel}
\usepackage{amssymb,amsmath}
\usepackage{cite}

\newcounter{par}
\newcounter{sbpr}
\newtheorem{thm}{Theorem}
\newtheorem{utt}{Proposition}
\newtheorem{cor}{Corollary}
\newtheorem{lm}{Lemma}
\theoremstyle{remark}
\newtheorem{rem}{Remark}
\theoremstyle{definition}
\newtheorem{df}{Definition}
\newtheorem{ex}{Example}

 \def\T{{\mathbb T}}
\def\C{{\mathbb C}}

\def\R{{\mathbb R}}

\def\P{{\mathbb P}}

\def\codim{{\rm codim\:}}

\def\1{{1}_\T}

\def\({{\rm (}}
\def\){{\rm )}}

\def\P{{\rm P}}
\def\res{{\rm res}}

\def\v{{\rm vol}}

\author{B. Kazarnovskii}
\title{On product of cocycles in a polyhedral complex\thanks{The research was carried out at the IITP RAS at the expense of the Russian Foundation for Sciences (project № 14-50-00150)}}
\begin{document}
\maketitle {\footnotesize
A product of cochains in a polyhedral complex is constructed.
The multiplication algorithm depends on the choice of a parameter.
The parameter is a linear functional on the ambient space.
Cocycles form a subring of the ring of cochains,
cobounds form an ideal of the ring of cocycles,
and the quotient ring is a ring of cohomologies.
If the complex is simplicial,
then the algorithm reduces to the well known algorithm of ${\rm\breve{C}}$ech.
Similar algorithms are used in geometry of polyhedra
for multiplying of cocycles taking values in the exterior algebra of the ambient space.
Therefore, we assume the ring of values to be supercommutative.
This text contains the basic definitions,
the statements of main theorems and some of their applications in convex geometry.
The full text will be offered to "Izvestiya: Mathematics".
}
\setcounter{par}{1}
\paragraph{\arabic{par}. Introduction.}
Let $X$ be a finite set of closed convex ($\leq k)$-dimensional
polyhedra in a real vector space $V$.
The polyhedra in $X$ are called \emph{the cells} of $X$.
We say that $X$ is a $k$-dimensional \emph{polyhedral complex} (\P-complex for short) if

(1) any face of any cell is also a cell

(2) any non-empty intersection of two cells is a common face of these cells.
\par\medskip
If any cell is a simplex then $X$ is called a simplicial complex.

Let $r$ be a function on the set of oriented cells
taking values in a supercommutative
ring $S$ (note that any commutative ring is also super\-com\-mutative).
We assume the function $r$ to be odd
(i.e. its value changes sign with the changing of the argument cell orientation).
Function $r$ is said to be a cochain or (depending on the context) a chain with values in $S$.
If the parity of all the values of cochain $s$ are the same,
then we denote the parity by $|s|$
(if $S$ is commutative, then $|s|=0$).
Let $C^p(X,S)$ or $C_p(X,S)$ be respectively an $S$-module of $p$-cochains or an
$S$-module of $p$-chains with values in $S$.
Define the boundary operators $d\colon C^p(X,S)\to C^{p+1}(X,S)$
and $\partial\colon C_p(X,S)\to C_{p-1}(X,S)$ as usual.
In what follows speaking about the cochains
(but not about the chains!) we assume
that any cell is compact.

For a general point $v$ of a dual space $V^*$
we define the product of cochains $r_p\smallsmile_vr_q\in C^{p+q}(X,S)$.
%For any general $u,v$ is true that $r_p\smallsmile_vr_q-r_p\smallsmile_ur_q\in dC^{p+q-1}(X,S)$.
For any choice of $v$ we get the ring of cochains.
The cocycles form a subring. % of the ring of cochains.
The cobounds form an ideal of the ring of cocycles.
The quotient of cocycles modulo the ideal of cobounds
does not depend on the choice of parameter and is isomorphic to the ring of cohomologies.

The use of superrings is motivated by geometric applications.
In these applications the cocycles take their values in the exterior algebra $\bigwedge^* V$.
It is important in geometry
that the product of some special cocycles
does not depend on the choice of the parameter.
Therefore, we act with the multiplication of cocycles
more careful,
than it is necessary for cohomologies.
Such an accuracy in the description of product
is also sometimes required in topology \cite{St}.

The parameter $v$ is controlled by the discriminant $D(X)$.
The discriminant is a union of subspaces of $V^*$
of codimension $1$.
The set of subspaces is finite.
If $k=1$ then $D(X)=\emptyset$.
If $k=2,3$ then the discriminant is a union of the subspaces
orthogonal to the edges ($1$-dimensional cells) of \P-complex $X$.
For all
$v$ in the connected component of $V^*\setminus D(X)$
the product of cocycles (not cochains!) $r_p\smallsmile_vr_q$ is the same.

Let $r_p\in C^p(X,S), r_q\in C^q(X,S)$.
Suppose first
that $X$ is a simplicial complex with ordered vertices
(the ordering of vertices is a parameter).
Eduard ${\rm\breve{C}}$ech in 1936 \cite{Ch}
proposed the following product of
cochains:
\begin{equation}\label{sympl}
 (r_p\smallsmile r_q)([u_0,\cdots,u_{p+q}])=r_p([u_0,\cdots,u_p])\:r_q([u_p,u_{p+1},\cdots,u_{p+q}]),
\end{equation}
where $[u_0,\cdots,u_m]$ is an $m$-dimensional simplex with the increasing vertices $u_0,\cdots,u_m$.
It is proved in \cite{Wh} that (\ref{sympl})
gives the (usual in topology) product of cohomologies.
We give an analogue of the ${\rm\breve{C}}$ech algorithm for product of cochains in an arbitrary \P-complex $X$.

Let $\gamma\in X$ and $\dim\gamma=p+q$. Then
\begin{equation}\label{formula}
 (r_p\smallsmile_vr_q)(\gamma) = \sum_{(\delta,\lambda)\in{\cal P}(p,q,\gamma,v)} r_p(\delta)r_q(\lambda),
\end{equation}
where ${\cal P}(p,q,\gamma,v)$ is some subset
of the set of pairs of oriented faces $(\delta,\lambda)$ of $\gamma$
of dimensions $p,q$.
An algorithm for the choose of the subset ${\cal P}(p,q,\gamma,v)$
is given in \S 2.
Assume
that the vertices of $X$ are numbered in ascending values
of functional $v$.
If
$\gamma$ is a simplex
then the set ${\cal P}(p,q,\gamma,v)$ consists of one element
and (\ref{formula}) turns to (\ref{sympl}).
The formulas similar to (\ref{formula}) are used in the geometry of polyhedra, in tropical geometry,
in the theory of toric varieties and in complex analysis
\cite{Ask1, Fult, K2, K, Ask3, Est, K3}.

In \S 2 we formulate the main theorems
(Theorem \ref{thm1}, \ref{thm2} and \ref{thm3})
and the necessary definitions.
The proofs of theorems are based on the structure of local duality.
The duality transforms the product of cocycles into the intersection of cycles.
For some of \P-complexes there exists \emph{a geometrical construction} of the dual \P-complex.
A cochain complex of the \P-complex is isomorphic to a chain complex of the dual \P-complex.
Any \P-complex is composed of \ P-subcomplexes
for which the dual \P-complexes exist.
Respectively a cochain complex of any \P-complex is glued of chain complexes
of locally dual \P-complexes.
So the proof of theorems
reduces to
the corresponding statements about the intersections of cycles of \P-complexes.
These statements are of interest in themselves.
On their application in geometry see \cite{K2, K, Est, K3}.

The proof of Theorems \ref{thm1} and \ref{thm2} uses the simplest version of duality:
a \P-complex
formed by the faces of a convex polyhedron is dual to
the fan of cones dual to the polyhedron.
Respectively
the cochain complex of a convex polyhedron is isomorphic to
the chain complex of a dual fan of cones.
The proof of Theorem \ref{thm3} uses some more delicate construction of duality.
It is connected to a notion of a regular partition of polyhedron \cite {Ask3}.

This publication has two direct sources.
The first is the construction of the intersection of tropical varieties
as the intersection of cycles with values in exterior algebra \cite{K2,K}.
The second source is the unpublished text of A. Khovanskii,
in which
the tropical varieties represented as cocycles
of regular polyhedra partitions.
Thus we get a dual interpretation,
in which the intersection of varieties becomes a product of cocycles.
I am grateful to A. Khovanskii for many useful discussions
and comments on his text.

In the remaining part of the Introduction we discuss the simplest geometrical applications of (\ref{formula}).

Let $U$ be a bounded open set in a $q$-dimensional oriented affine subspace of the space $V$.
Let $\beta_{U}\in\bigwedge^qV$ be a multivector for which
$\int_{U}\omega=\omega(\beta_{U})$ for any constant exterior form $\omega\in\bigwedge^q V^*$.
Multivector $\beta_{U}$ changes sign with the changing of the orientation of $U$.
We view $\beta_{U}$ is a $q$-dimensional volume of $U$.

Let $X$ be a \P-complex in $V$ and $S=\bigwedge^*V$.
Define a cochain $\v^X_q\in C^q(X,S)$ as $\v^X_q(\delta)=q!\beta_{\delta}$.
The cochain $\v^X_q$ is a cocycle.
It is equivalent to the Pascal equations for the $(q+1)$-dimensional
cells of $X$.
(The Pascal equations for a k-dimensional polyhedron
 are of the form $\sum_\lambda v_\delta = 0$,
 where $\delta$ is $(k-1)$-dimensional face of the polyhedron
 and $v_\delta$ is the exterior normal
vector to $\delta$ of the length equal to the $(k-1)$-dimensional volume of $\delta$.)
 \begin{utt} \label{uttK}
 $\forall v\in V\colon\,\,\, \v^X_p\smallsmile_v\v^X_q=\v^X_{p+q}$.
 \end{utt}
 Let $V_\delta$ be a subspace of $V$ generated by the differences of points
of the cell $\delta$.
 Let ${\cal R}^p(X,S)$ be a set of cocycles $r_p\in C^p(X,S)$ for which $r_p(\delta)\in\bigwedge^pV_\delta$.
 For example $\v^X_q\in{\cal R}^q(X,S)$.
The product $\smallsmile_v$ defines the structure of commutative ring
on ${\cal R}^*(X,S)={\cal R}^0(X,S)\oplus\cdots\oplus{\cal R}^n(X,S)$
(commutativity follows from the Proposition \ref{utt1} (1) in \S 2).
\begin{utt} \label{uttTrop}
The ring ${\cal R}^*(X,S)$ does not depend on a choice of $v$.
 \end{utt}
If the \P-complex $X$ consists of the faces of a convex $n$-dimensional polyhedron,
then Proposition \ref{uttK} gives an algorithm for computing
of volumes and mixed volumes of polyhedra.
It follows from Proposition \ref{uttK}
that $\v^X_n=(\v^X_1)^n$.
Therefore (\ref {formula}) leads to an algorithm for calculating the volume of a polyhedron with
a given set of edges \cite{Ask1}. The algorithm consists of

(1) identification of some sets consisting of $n$ edges of the polyhedron
(sets of edges depend on the choice of the parameter $v\in V^*$)

(2) summation of volumes of the corresponding parallelepipeds.

A similar algorithm to compute the mixed volume is as follows.
Let $\gamma$ be a Minkowski sum of polyhedra $\gamma_1,\cdots,\gamma_{n}$.
We define the $1$-cocycles $\v^X_{1,k}$,
setting their values on the edges of $\gamma$ by the following rule.
Any edge $\lambda\subset\gamma$ is uniquely represented as a Minkowski sum
$\lambda_1+\cdots+\lambda_n$,
where $\lambda_k$ is either an edge or a vertex
of $\gamma_k$.
If $\lambda_k$ is the vertex then set $\v^X_{1, k}(\lambda) = 0$,
if $\lambda_k$ is the edge then set $\v^X_{1, k}(\lambda) = \lambda_k$.
%Let $\v^X_{1, k}(\lambda) = 0$,
%if $\lambda_k$ is the vertex,
%and let $\v^X_{1, k}(\lambda) = \lambda_k$,
%if $\lambda_k$ is the edge.
Then the $n$-vector $\v^ X_ {1,1}\smallsmile\cdots\smallsmile\v^X_{1, n}$
is a mixed volume of polyhedra $\gamma_1,\cdots,\gamma_{n}$,
multiplied by $n!$.
\begin{utt} \label{uttToric}
Let $\gamma$ be an integer simple convex polyhedron and let $M$
be a corresponding toric variety {\rm\cite{Ask2}}.
Then the ring ${\cal R}^*(X, S)$ coincides with
the cohomology ring of $M$ {\rm\cite{Fult}}.
 \end{utt}
Let $X$ be a \P-complex in $V$ and
$a_1,\cdots,a_p$ are $0$-cochains with values in $V$.
Suppose that $da_i\in{\cal R}^1(X,S)$.
For example,
if $a(\lambda)=\lambda$,
then $da=\v^X_1\in{\cal R}^1(X,S)$.
Consider the product $a_1\cdots a_p$
as a $0$-cochain $\prod(a_1,\cdots,a_p)$ with values in symmetric algebra of the space $V$.
\begin{utt} \label{uttEst}
The product of cocycles $da_1\smallsmile\cdots\smallsmile da_p\in C^p(X,\bigwedge^* V)$
%with values in $\bigwedge^* V$
is completely determined by the $0$-cochain $\prod(a_1,\cdots,a_p)$.
\end{utt}
Proposition \ref{uttEst} 
reduces 
(by (\ref{formula}) and Proposition \ref{uttTrop})
to the case
of \P-complex $X$
consisting of the faces of a convex polyhedron.
In this case the above algorithm for computing of the mixed volume together with Proposition \ref{uttEst}
gives the following corollary.
\begin{utt} \label{uttEst2}
A mixed volume of polyhedra is completely determined
by the product of the support functions of polyhedra {\rm\cite{Est}}.
\end{utt}
The last statement is not true for arbitrary convex bodies \cite{K3}.

Let $X$ be a \P-complex in $\C^n$.
In this case
the multiplication of cocycles with values in $\bigwedge^*_\C\C^n$
also leads to some geometric applications.
For example let $\C^n_\R$ be a realification of $\C^n$ and
let $\rho\colon\bigwedge^*_\R\C^n_\R\to\bigwedge^*_\C\C^n$ be a
unique homomorphism of rings,
such that $\rho\colon\C^n_\R\to\C^n$ is the
identity map.
Define a cocycle $\v^{X,\C}_p$ with values in $\bigwedge^*_\C\C^n$
as $\v^{X,\C}_p(\gamma) =\rho(\v^X_p(\gamma))$.
Propositions \ref{uttK}, \ref{uttEst}
for cocycles $\v^{X,\C}_p$
also are true.
Their consequences are some of the properties of the mixed pseudovolume
of polyhedra in a complex space \cite{Onul, ExpNewt, Al}.
These properties similar to the above properties of the mixed volume.
The products of cocycles $\v^{X,\C}_p$
is also related to the action of the complex Monge-Ampere operator
on piecewise linear functions in the space $\C^n$ \cite{K, K3}.
\addtocounter{par}{1}
\paragraph{\arabic{par}. Basic definitions and main theorems.}
The construction of the set ${\cal P}(p,q,\gamma,v)$ from (\ref{formula})
uses the definition of a fan of cones
dual to the convex polyhedron.
(\P-complex is called a fan of cones
if its cells are convex cones with the vertex $0$).
Let $V_\delta$ be a subspace of $V$ generated by the differences of points
of the cell $\delta$.

Let $\gamma$ be a convex polyhedron in $V$,
$\dim V=n$.
By definition the dual cone $\delta^*$ of a face $\delta$
consists of all points $u\in V^*$ such that $\max_{v\in\gamma}u(v)$
is attained for any $v\in\delta$.
The cone $\delta^*$
belongs to the orthogonal complement of the subspace $V_\delta$ and
$\dim\delta^*=n-\dim\delta$.
\begin{rem}[local duality]\label{remLocDualSimple}
Let $Y$ be a \P-complex consisting of all faces of $\gamma$ and let
$Y^*$ be a fan of cones dual to $\gamma$.
If the space $V$ is oriented
then the orientations of the cells $\delta, \delta^*$ can be agreed.
For $r\in C^p(Y,S)$ define the chain $r^*\in C_{n-p}(Y^*,S)$
as $r^*(\delta^*)=r(\delta)$.
The mapping $r\mapsto r^*$ is an isomorphism of complexes
$C^*(Y, S)\to C_*(Y^*,S)$.
\end{rem}
\begin{df}\label{dfIntersP}
For \P-comlexes $\mathfrak{X}_1,\cdots,\mathfrak{X}_k$ define the intersection $\mathfrak{X}_1\cap\cdots\cap\mathfrak{X}_k$ as
a \P-complex consisting of the cells $\delta_1\cap\cdots\cap\delta_k$ with $\delta_i\in\mathfrak{X}_i$.
\end{df}
\begin{df}\label{dfTransv}
We say that \P-complexes $\mathfrak{X}_1,\cdots,\mathfrak{X}_k$ in the space $V$
are trans\-versal,
if for any cells $\delta_i\in\mathfrak{X}_i$ with non-empty intersection
$\delta_1\cap\cdots\cap\delta_k$
it is true that $\codim\cap_iV_{\delta_i}=\sum_i\codim V_{\delta_i}$.
\end{df}
Let $X$ be a \P-complex in $V$.
Let $\gamma\in X$ and let $X_\gamma$ be a \P-complex in $V_\gamma$,
consisting of the faces of $\gamma$.
Let $\dim\gamma=p+q$.
We denote by $(X^*_\gamma)^q$ and $(X^*_\gamma)^p$, respectively, $q$-dimensional and $p$-dimensional
skeletons of the fan of cones $X^*_\gamma\subset(V_\gamma)^*$ dual to the polyhedron $\gamma$.
\begin{df}[convenient  functional]\label{dfGoodPoint}
A point $v\in V^*$ is said to be $(p,q)$-\emph{convenient point},
if \P-complexes $\pi_\gamma(v) + (X^*_\gamma)^q$ and $(X^*_\gamma)^p$ are transversal
for any $\gamma\in X$
where $\pi_\gamma\colon V^*\to(V_\gamma)^*$ is a projection
conjugate to the imbedding $V_\gamma\to V$.
Let ${\cal U}(X)=\cap_{p,q}{\cal U}_{p,q}(X)$,
where ${\cal U}_{p,q}(X)$ is a set of $(p,q)$-convenient points.
The points of ${\cal U}(X)$ are said to be convenient.
\end{df}
\textbf{The algorithm for the selecting of the subset ${\cal P}(p,q,\gamma,v)$ in formula} (\ref{formula}) is as follows.
If $v\in{\cal U}_{p,q}(X)$
 then the \P-complex $\left(\pi_\gamma(v)+(X^*_\gamma)^q\right)\cap (X^*_\gamma)^p$
 %является конечным множеством
 is a finite set of points.
 It consists of the points $(\pi_\gamma(v)+\delta^*)\cap\lambda^*$
 where $(\delta,\lambda)$ is a pair of faces of the polyhedron $\gamma$
 of dimensions $(p,q)$
 (the set of pairs depends on the choice of $v$).
 Choose the orientations of the cells $\delta,\lambda$ so
 that the orientations of the spaces $V_\gamma$ and $V_\delta\oplus V_\lambda$
 agreed with the isomorphism of vector spaces $V_\gamma = V_\delta\oplus V_\lambda$.
 Let ${\cal P}(p,q,\gamma,v)$ be a resulting set of ordered pairs of oriented cells.
 %The terms of the formula (\ref{formula}) are defined.
 \emph{If the point $v$ is convenient,
 then the set ${\cal P}(p,q,\gamma,v)$ and the $v$-product $r\smallsmile_vs$ of any cochains $r,s$ constructed}.
\begin{lm} \label{lmC}
If $\gamma$ is a simplex
then the set ${\cal P}(p, q, \gamma, v)$ consists of one element $(\delta,\lambda)$,
$\dim\delta\cap\lambda = 0$
and $v(x)\leq v(y)$ for all $x\in\delta, y\in\lambda$.
\end{lm}
\begin{cor}\label{corC}
Assume
that the \P-complex $X$ is simplicial and its vertices are numbered in ascending values
of functional $v$.
Then
the products of cochains
by the formulas {\rm(\ref {sympl})} and {\rm(\ref {formula})} coincide.
\end{cor}
\begin{cor}\label{corAfterDef}
Let the points $u, v$ belong to the same connected component
of ${\cal U}_{p, q}(X)$.
Then $r_p\smallsmile_ur_q = r_p\smallsmile_vr_q$
for any $r_p\in C^p(X, S)$, $r_q\in C^q(X, S)$.
\end{cor}
\begin{utt}\label{utt1}
Let $v\in{\cal U}_{p,q}(X)$, $r_p\in C^p(X,S)$, $r_q\in C^q(X,S)$.
Then

{\rm(1)} $r_p\smallsmile_vr_q=(-1)^{pq+|r_p||r_q|}r_q\smallsmile_vr_p$.

{\rm(2)} $d(r_p\smallsmile_vr_q)=dr_p\smallsmile_vr_q + (-1)^{p} r_p\smallsmile_vdr_q$.
 \end{utt}
\begin{cor}\label{cor1}
Let the cochains $ r_p, r_q $ are closed. Then

{\rm(1)} the cochain $r_p\smallsmile_vr_q$ is closed;

{\rm(2)} if $r_p=dr_{p-1}$, then $r_p\smallsmile_vr_q=d(r_{p-1}\smallsmile_vr_q)$.
 \end{cor}
\begin{df}[special triples of cells\footnote
{close to the notion of "$1$-regular pair of simplexes"\ in \cite{St}.}]\label{dfPairSpec}
We denote by $\Lambda_{p, q}$ a set of triples of cells $(\delta,\lambda,\mu)$
of dimensions $p,q,p+q$ such
that $\delta\cup\lambda\subset\mu$ and
$\dim V_\delta\cap V_\lambda\geq1$.
If $(\delta,\lambda,\mu)\in\Lambda_{p,q}$,
$\dim V_\delta\cap V_\lambda=1$ and $\delta\cup\lambda\subset\gamma$,
where
$\gamma$ is a $(p+q-1)$-dimensional face of the polyhedron $\mu$,
then the triple
$(\delta,\lambda,\gamma)$ is said to be \emph{a special triple of cells}.
\end{df}
\begin{ex}\label{exSpec}
Let $\gamma\in X$ and $\gamma$ is not a maximal cell of \P-complex $X$.
Then the following triples $(\delta,\lambda,\gamma)$ are special\\
$\dim\gamma=1\colon$ $(\gamma,\gamma,\gamma)$;\\
$\dim\gamma=2\colon$ $(\delta,\gamma,\gamma)$
 where $\delta$ is a side of the polygon $\gamma$;\\
$\dim\gamma=3\colon$ $(\delta,\gamma,\gamma)$,
 where $\delta$ is an edge,
 and $(\delta,\lambda,\gamma)$,
 where $\delta,\lambda$ is a pair of not parallel $2$-dimensional faces of $\gamma$.
\end{ex}
\begin{df}[discriminant of \P-complex]\label{dfDiscr}
Let $V^*_{\delta,\lambda}$ be a subspace of $V^*$ orthogonal to $V_\delta\cap V_\lambda$.
The union of all hyperplanes $V^*_{\delta,\lambda}$,
taken over all the special triples $(\delta,\lambda,\gamma)$,
denote by $D(X)$ and call \emph{a discriminant} of $X$.
\end{df}
 \begin{thm}\label{thm1}
 Let the cochains $ r_p, r_q $ are closed,
 the points $u,v$ are convenient and $u,v\in{\cal C}$,
 where ${\cal C}$ is a connected component of $V^*\setminus D(X)$.
Then $r_p\smallsmile_ur_q=r_p\smallsmile_vr_q$.
\end{thm}
\begin{df}\label{dfBadPlanes}
Let $(\delta,\lambda,\mu)\in\Lambda_{p,q}$.
If $\dim V_\delta\cap V_\lambda=1$
then the hyperplane $V^*_{\delta,\lambda}\subset V^*$
is said to be \emph{uncovenient}.
Uncovenient hyperplane $H$ is said to be \emph{the discriminant hyperplane}
if $H\subset D(X)$.
If $\dim V_\delta\cap V_\lambda>1$ and $(u+\delta^*)\cap\lambda^*\not=\emptyset$,
where $\delta^*,\lambda^*$ are dual to $\delta,\lambda$,
then the point $u$ is said to be \emph{uncommon}.
\end{df}
\begin{cor}\label{corUnCommon}
A codimension of the set of uncommon points $>1$.
\end{cor}
\begin{cor}\label{corBadPlanes}
Let $D'(X)$ be a set of all not convenient points,
which are not uncommon.
Then

{\rm(1)} $D'(X)$ is an open set in the union of all unconvenient hyperplanes

{\rm(2)}$D(X)\subset D'(X)$

{\rm(3)} $\codim D'(X)\setminus D(X)=1$.
\end{cor}
Let $\delta,\lambda,\gamma$ be a special triple and $\phi\not\in V^*_{\delta,\lambda}$.
Define an agreement for orien\-tations of the cells $\delta,\lambda,\gamma$ as follows.
The functional $\phi$ gives the orientation of $1$-dimensional subspace $V_\delta\cap V_\lambda$.
Choose orientations of the the cells so
that the isomorphism
$V_\gamma/(V_\delta\cap V_\lambda)=V_\delta/(V_\delta\cap V_\lambda)\oplus V_\lambda/(V_\delta\cap V_\lambda)$
is an isomorphism of oriented spaces.
\begin{df}[special cochains]\label{dfCochainSpec}
Let $r_p\in C^p(X,S)$, $r_q\in C^q(X,S)$ and let $(\delta,\lambda,\gamma)$ be a special triple,
$\dim\delta=p$, $\dim\lambda=q$ and $\phi\not\in V^*_{\delta,\lambda}$.
We denote by $\vartheta^{\delta,\lambda,\gamma;\phi}_{r_p,r_q}$ $\,(p+q-1)$-cochain
equal to $r_p(\delta) r_q(\lambda)$ on the cell $\gamma$
(here we use the above-described agreement of orientations that uses the choice of $\phi$)
and equal to $0$ at all other cells of $X$.
\end{df}
Let $u,v,\kappa\in V$ and let it is true that

(1) $([u,v]\setminus\kappa)\subset{\cal U}_{p,q}$,
where $\kappa$ is the inner point of a segment $[u,v]$

(2) $\kappa\in H$, where $H$ is a discriminant hyperplane

(3) $\kappa$ is not uncommon

(4) if $\kappa\in G$,
where $G$ is an uncovenient hyperplane, then $G=H$.

If
$H=V^*_{\delta,\lambda}$,
then (4) implies
that $(\kappa+\delta^*)\cap\lambda^*\subset H$,
where $\delta^*,\lambda^*$ are dual to the faces $\delta,\lambda$ of a polyhedron $\gamma$.
It is also true that if $(\kappa+\delta^*)\cap\lambda^*\not=\emptyset$,
then $H=V^*_{\delta,\lambda}$.
From (2-3) it follows
that the intersection $(\kappa+\delta^*)\cap\lambda^*$ is either empty or
transversel in the space $H$.
The set of special triples $(\delta,\lambda,\gamma)$ such
that $(\kappa+\delta^*)\cap\lambda^*\not=\emptyset$
and $\dim\delta=p$, $\dim\lambda=q$,
denote
${\cal S}^H_{p, q}(\kappa)$.
\begin{thm}\label{thm2}
Let
$r_p,r_q$ are the cocycles of dimensions $(p,q)$,
${\cal C},{\cal D}$ are the
connected components of $V^*\setminus D(X)$
neighboring on $(n-1)$-dimensional cone of the subspace $H$,
$u\in{\cal C}, v\in{\cal D}$.
Then
\begin{equation}\label{cocycle1}
 r_p\smallsmile_{v}r_q - r_p\smallsmile_{u}r_q=\sum_{(\delta,\lambda,\gamma)\in{\cal S}^{H}_{p,q}(\kappa)} d\vartheta^{\delta,\lambda,\gamma;v-u}_{r_p,r_q}.
\end{equation}
\end{thm}
\begin{cor}\label{corMain2}
For any convenient points $u,v$ and for any cocycles $r,s$ it is true that
$r\smallsmile_us - r\smallsmile_vs\in dC^*(X,S)$.
\end{cor}
\begin{cor}\label{corMain}
For any convenient point $v$ the product $\smallsmile_v$ defines
the product of cohomologies
independent of the choice of $v$.
\end{cor}
\begin{rem}\label{remDiscrCone}
Theorem \ref{thm1} implies that the cochain
$\sum_{(\delta,\lambda,\gamma)\in{\cal S}^{H}_{p,q}(\kappa)} d\vartheta^{\delta,\lambda,\gamma;v-u}_{r_p,r_q}$
from (\ref{cocycle1})
(as opposed to the set ${\cal S}^H_{p, q}(\kappa)$)
does not depend on the choice of $\kappa$.
\end{rem}
\begin{rem}\label{remIndepend}
Let $S=\bigwedge^*V$ and $r_p\in{\cal R}^p(X,S)$, $r_q\in{\cal R}^q(X,S)$
(see \S 1).
Then for any special triple $(\delta,\lambda,\gamma)$ it is true that
$r_p(\delta)r_q(\lambda)=0$ and hence $\vartheta^{\delta,\lambda,\gamma;\phi}_{r_p,r_q}=0$.
This implies that the Proposition \ref{uttTrop} follows from theorems \ref{thm1} and \ref{thm2}.
\end{rem}
\begin{df}\label{dfRes}
Let \P-complex $X$ be a subdivision of \P-complex $Y$.
By definition the map $\res\colon C^*(X,S)\to C^*(Y,S)$
acts on a cochain $r_p\in C^p(X,S)$ as
$\res(r_p)(\gamma)=\sum_{X\ni\delta\subset\gamma} r_p (\delta)$
where $\gamma\in Y$ is an arbitrary $p$-dimensional cell.
\end{df}
\begin{cor}\label{corCoBeta}
The map $\res$ is a homomorphism of cochain complexes.
\end{cor}
\begin{utt}\label{uttRes}
Let the \P-complex $X$ is simplicial.
Then the map $\res\colon C^*(X,S)\to C^*(Y,S)$
induces an isomorphism of cohomologies $H^*(X,S)\to H^*(Y,S)$.
\end{utt}
\begin{thm}\label{thm3}
Let the \P-complex $X$ is simplicial,
$v\in{\cal U}(Y)\cap{\cal U}(X)$ and $r_p,r_q$ are the cocycles of
%\P-complex
$X$.
Then $\res(r_p)\smallsmile_v\res(r_q) - \res(r_p\smallsmile_v r_q)\in dC^{p+q-1}(Y,S)$.
\end{thm}
It follows from Theorem \ref{thm3} and Proposition \ref{uttRes}
that $\smallsmile_v$-product of cohomo\-logies with any convenient $v$
coincides with the Kolmogorov-Alexander product.
\begin{rem}\label{remPartition}
The proof of Theorem \ref{thm3} implies
that if $S=\bigwedge^* V$ and $r_p,r_q\in{\cal R}^*(X, S)$
(see Proposition \ref{uttTrop}),
then $\res(r_p)\smallsmile\res(r_q)=\res(r_p\smallsmile r_q)$.
For cocycles $\v^Y_p$ (see \S 1) it is true that $\res(\v^X_p)=\v^Y_p$.
Therefore Proposition \ref{uttK} is reduced to
the \P-complex $X$ consisting of the faces of a simplex.
\end{rem}
\begin{rem}\label{remToThm3}
For the cobound $\res(r_p)\smallsmile_v\res(r_q)-\res(r_p\smallsmile_vr_q)$ there exists a formula similar to
(\ref {cocycle1}) from Theorem \ref{thm2}.
\end{rem}
%
%
%
%
%\vfill\eject
\selectlanguage{english}

\par\medskip
\begin{flushleft}
{\small
Institute for Information Transmission Problems, Moscow, Russia\\
E-mail address: kazbori@iitp.ru
}
\end{flushleft}

\end{document}